\documentclass [11pt] {article}

\usepackage {amssymb}

\title {A Remark on Non-equivalent Star Products \\
        via Reduction for $\CP$}

\author {{\bf Stefan 
          Waldmann\thanks{Stefan.Waldmann@physik.uni-freiburg.de}
         } \\[3mm]
         Fakult\"at f\"ur Physik\\Universit\"at Freiburg \\
         Hermann-Herder-Str. 3 \\
         79104 Freiburg i.~Br., F.~R.~G \\[3mm]
        }

\date{FR-THEP-98/3 \\[1mm]
      17 February 1998 \\[5mm]}

\newcommand {\C} {{\mathbb C^{n+1}\setminus\{0\}}}
\newcommand {\CP} {{\mathbb {CP}^n}}
\newcommand {\wt} [1] {{\widetilde {{#1}}}}

\newcommand {\im} {{\bf i}}
                                  
\newcommand {\BEQ} [1] {\begin {equation} \label {#1}}
\newcommand {\EEQ} {\end {equation}}

\newcommand {\cc} [1] {\overline {{#1}}}

\newcommand {\ad} {{\rm ad}}

\newenvironment {proof}{\small {\sc Proof:}}{{\hspace*{\fill} $\square$}}

\newtheorem {lemma} {Lemma} [section]

\newtheorem {theorem} [lemma] {Theorem}
\newtheorem {corollary} [lemma] {Corollary}

\begin {document}

\maketitle

\begin {abstract}
In this paper we construct non-equivalent star products on $\CP$ by phase 
space reduction. It turns out that the non-equivalent star products occur 
very natural in the context of phase space reduction 
by deforming the momentum map of the $U(1)$-action on $\C$ 
into a quantum momentum map and the corresponding momentum value into a 
quantum momentum value such that the level set, 
i.~e. the `constraint surface', of the quantum momentum map coincides 
with the classical one. All equivalence classes of star products on $\CP$
are obtained by this construction.
\end {abstract}


\section {Introduction}
\label {IntroSec}

The concept of deformation quantization as introduced by Bayen, Flato, 
Fr{\o}nsdal, Lichnerowicz, and Sternheimer in \cite {BFFLS78} is now a
well-established way to understand quantization of classical systems: the 
algebra of classical observables, i.~e. the smooth complex-valued
functions $C^\infty (M)$ on a symplectic manifold, the phase space of the 
system is deformed by introducing an associative formal so-called star 
product $*$ for $C^\infty (M)[[\lambda]]$ depending on a formal parameter 
$\lambda$ such that the zeroth order of the star product is the pointwise 
product and the $*$-commutator of two functions equals in first order 
$\im$ times the Poisson bracket. Hence the formal parameter $\lambda$ 
is to be identified with Planck's constant $\hbar$ and the algebra of 
quantum observables $C^\infty (M)[[\lambda]]$ turns out to be a 
deformation of the classical one in the sense of Gerstenhaber 
\cite {GS88}. The existence of such star products for symplectic
manifolds was shown by DeWilde and Lecomte \cite {DL83},
Fedosov \cite {Fed94}, and Omori, Maeda, and Yoshioka \cite {OMY91} and
recently the existence
of star products for arbitrary Poisson manifolds was stated by Kontsevich 
\cite {Kon97}. The classification of star products up to equivalence by 
formal power series with coefficients in the second de Rham cohomology 
was shown by Nest and Tsygan \cite {NT95a} and by Bertelson, Cahen, and Gutt 
\cite {BCG97}.

In this paper we continue the work done together with Bordemann, Brischle 
and Emmrich in \cite {BBEW96a,BBEW96b} where the star product analogue of the 
Marsden-Weinstein phase space reduction for the example of the 
$U(1)$-reduction of $\C$ to $\CP$ was considered. Apart form a few papers 
the subject of reduction for star product seems not to be studied very 
intensively until now: In \cite {BFFLS78} several basic examples were 
considered and in \cite {BBEW96a,BBEW96b} the example of $\CP$ and its 
non-compact dual was considered and the method of reduction was applied 
to find explicit formulas for star products on $\CP$. A generalization to 
complex Grassmannian manifolds is given in \cite {Sch97} by Schirmer. 
Recently Fedosov gave a general construction for the reduction in case 
of a Hamiltonian group action of an arbitrary compact Lie group in 
\cite {Fed97}. All these reductions proceed more or less the same: one 
starts with a suitable star product on the `big phase space' such that the 
invariant functions form a subalgebra. Then the ideal generated by the 
components of the momentum map minus the momentum value is factored out 
and it remains to show that the quotient algebra is isomorphic to the 
functions on the reduced phase space endowed with a suitable star product 
which is hence called the reduced star product. This construction is 
physically reasonable and provides a way to perform a reduction in 
the quantum case too.

Nevertheless none of these approaches seems to deal with the question 
whether an equivalence transformation in the big phase space results in an 
equivalent reduced star product. In our example the possibility of 
non-equivalent star products has to be taken into account since the second 
de Rham cohomology of $\CP$ is one-dimensional. This work, originally 
motivated by a discussion with Bordemann, Flato, Schirmer and Sternheimer,
provides a very simple example that non-equivalent star products may 
occur. It depends crucially on the definition of the above mentioned
ideal and here one has in principle at least two reasonable possibilities:
On the one hand for a fixed star product on the big phase space
one can fix the value of the corresponding quantum momentum map to 
be the classical one which results in our example in only one star 
product for the quotient independent of the chosen star product on 
the big phase space. On the other hand one can fix the level set, 
i.~e. the `constraint surface', of the classical momentum map and hence one 
eventually has to modify the value of the (quantum) momentum map by 
`quantum corrections'. In this case it turns out that the resulting reduced 
star products are no longer equivalent in our example. Now both possibilities 
have their physical motivation and hence this example shows that there may 
occur some subtilities depending on the point of view which structure is
more important: the value of the classical momentum map or the classical 
constraint surface. Moreover this construction shows that non-equivalent 
star products arise very natural in the reduction process.

The paper is organized as follows: in section \ref {PreSec} we remember 
briefly the notation and results from \cite {BBEW96a}.
In section \ref {InvSec} we describe the invariance properties and the 
reduction of the star products and in section \ref {NonSec} we prove the 
non-equivalence of the star products obtained for $\CP$. 

\section {Preliminary results}
\label {PreSec}

Let us first remember the construction of star products of Wick type by 
phase space reduction for $\CP$ where we mainly use the notion as in 
\cite {BBEW96a}. We start with the K\"{a}hler manifold $\C$ 
with the usual K\"{a}hler form $\omega = \frac{\im}{2} dz^k \wedge d\cc z^k$
where $z^0, \ldots, z^n$ are the canonical holomorphic coordinates for 
$\C$ and summation over repeated indices is understood.
Moreover we consider the complex projective space $\CP$ and denote by 
$\pi: \C \to \CP$ the canonical (holomorphic) projection which maps 
$z \in \C$ to the ray $\pi (z) \in \CP$ through $z$. On $\C$ one has 
the usual $U(1)$-action $(e^{\im\varphi}, z) \mapsto e^{\im\varphi}z$ and 
the $\mathbb C\setminus\{0\}$-action $(\alpha, z) \mapsto \alpha z$. 
A function $f \in C^\infty (\C)$ is called {\em homogeneous} iff it is 
invariant under the $\mathbb C\setminus\{0\}$-action which is the case 
iff there exists a function $\phi \in C^\infty (\CP)$ such that 
$f = \pi^* \phi$. Moreover we consider the function 
$x : \C \to \mathbb R^+$ defined by $x (z) := \cc z^k z^k$.
Then a function $R \in C^\infty (\C)$ is called {\em radial} iff 
there exists a function $\varrho \in C^\infty (\mathbb R^+)$ such 
that $R = \varrho \circ x$.

Let us now recall the classical Marsden-Weinstein phase space reduction 
procedure for $\CP$ as e.~g. in \cite [p. 302] {AM85} to establish 
our notation: the $U(1)$-action is generated by the $\ad^*$-equivariant 
momentum map 
$J := -\frac{1}{2} x$ and any $\mu \in \mathbb R^-$ is a regular
value of $J$. Moreover $J^{-1} (\{\mu\})$ is just the $2n+1$ sphere 
centered at the origin with radius $\sqrt{-2\mu}$ in $\C$. Fix now once and 
for all an arbitrary value $\mu \in \mathbb R^-$ and denote by 
$i_\mu : J^{-1} (\{\mu\}) \to \C$ the inclusion and by 
$\pi_\mu: J^{-1} (\{\mu\}) \to J^{-1} (\{\mu\}) / U(1) \cong \CP$ the 
projection onto the reduced phase space.
Then $i^*_\mu \omega = \pi^*_\mu \omega_\mu$ determines the reduced
symplectic form $\omega_\mu$ on $\CP$
and it turns out that $\omega_\mu$ is just a multiple (depending on $\mu$) 
of the usual Fubini-Study form. For a $U(1)$-invariant function
$F \in C^\infty (\C)$ one defines the reduced function $F_\mu$ by
$F_\mu ([z]) := F \circ i_\mu (z)$ where $z \in J^{-1} (\{\mu\})$.

A suitable starting point for the deformation quantization of this
reduction is the {\em Wick star product} on $\C$ which is given by
the formal power series in $\lambda$ for $F, G \in C^\infty (\C)$ by
\BEQ {WickDef}
    F * G := \sum_{r=0}^\infty \frac{\lambda^r}{r!} 
             \frac{\partial^r F}{\partial z^{i_1} \cdots \partial z^{i_r}}
             \frac{\partial^r G}
             {\partial \cc z^{i_1} \cdots \partial \cc z^{i_r}}
\EEQ
which is known to be an associative formal star product for $\C$. For a 
more general treatment of this kind of star products on arbitrary K\"{a}hler 
manifolds see e.~g. \cite {BW97}. Note that in this normalization the
formal parameter $\lambda$ corresponds to $2\hbar$. Moreover we define
$\mathcal A := C^\infty (\C)[[\lambda]]$. Remember also the definition
of the bidifferential operators $M_r$ and $\wt M_r$ as introduced in
\cite [Eqn. 5 \& 23]{BBEW96a}: For $F, G \in C^\infty (\C)$ one defines
\BEQ {MrDef}
    M_r (F, G) := x^r 
    \frac{\partial^r F}{\partial z^{i_1} \cdots \partial z^{i_r}}
    \frac{\partial^r G}{\partial \cc z^{i_1} \cdots \partial \cc z^{i_r}}
\EEQ
and since for $\phi, \psi \in C^\infty (\CP)$ clearly 
$M_r (\pi^* \phi, \pi^* \psi)$ is again a homogeneous function
$\pi^* \wt M_r (\phi, \psi) = M_r (\pi^*\phi, \pi^* \psi)$
uniquely defines a bidifferential operator $\wt M_r$ on $\CP$. Crucial for 
the following is the observation that the formal power series
with coefficients in the $U(1)$-invariant functions on $\C$ which we shall 
denote by $\mathcal A^0 \subset \mathcal A$ build a sub-algebra 
with respect to the Wick product.

Let now
\BEQ {DDef}
    D (\lambda) := 1 + \sum_{r=1}^\infty \lambda^r d_r ,
    \quad
    C(\lambda) := D(\lambda)^{-1} = 1 + \sum_{r=1}^\infty \lambda^r c_r ,
    \quad
    d_r, c_r \in \mathbb C
\EEQ
be an arbitrary complex formal power series starting with $1$ and denote by 
$C(\lambda)$ the inverse series. Then for any such $D$ 
a formal series of differential operators $S_D : \mathcal A \to \mathcal A$ 
was constructed in \cite [Theorem 3.1] {BBEW96a} having the following 
properties: $S_D$ acts trivial on the homogeneous functions and 
$S_D : \mathcal A^0 \to \mathcal A^0$ and
\BEQ {SxDx}
    S_D x = D\left(\frac{\lambda}{x}\right) x .
\EEQ
This operator was used to define an equivalent star product which 
we shall now denote by $*^D$ to emphasise the dependence on $D$ defined by
\BEQ {StarDDef}
    F *^D G := S_D \left( (S_D^{-1} F) * (S_D^{-1} G)\right)
\EEQ
for $F, G \in \mathcal A$. It follows that the $*^D$-product of a radial 
function with an arbitrary $U(1)$-invariant function is just the pointwise 
product and that for homogeneous functions $f, g \in \mathcal A^0$ the 
equation 
\BEQ {HomProd}
    f *^D g = \sum_{r=0}^\infty \frac{1}{r!} 
    \left( \frac{\lambda}{x} C\left(\frac{\lambda}{x}\right)\right)^r
    \prod_{s=0}^r \left( 1 + s \frac{\lambda}{x} 
    C \left(\frac{\lambda}{x}\right)\right)^{-1} M_r (f, g)
\EEQ    
holds \cite [Eqn. 4]{BBEW96a}. Obviously this can be rearranged such that 
\BEQ {KrDef}
    f *^D g = \sum_{r=0}^\infty \left(\frac{\lambda}{x}\right)^r 
    K^D_r (f, g)
\EEQ
with some bidifferential operators $K_r^D$ which are linear 
combinations of the $M_r$ depending on the choice of $D$. Moreover 
$K_r^D (f, g)$ is clearly again homogeneous for $f, g$ homogeneous and 
thus there are again uniquely determined bidifferential operators 
$\wt K_r^D$ on $\CP$ such that
$\pi^* \wt K_r^D (\phi, \psi) = K_r^D (\pi^*\phi, \pi^*\psi)$
for $\phi, \psi \in C^\infty (\CP)$. By a simple computation using
the associativity of $*^D$ and the fact that the $*^D$ product of the 
radial functions $x^{-r}$ with any $U(1)$-invariant function is only the 
pointwise product one obtains that for $\phi, \psi \in C^\infty (\CP)$
\BEQ {RedStarDef}
    \phi *^D_\mu \psi := \sum_{r=0}^\infty 
    \left(\frac{\lambda}{-2\mu}\right) \wt K_r^D (\phi, \psi)
\EEQ
is an associative star product for $(\CP, \omega_\mu)$ and clearly 
$\left( \pi^*\phi *^D \pi^* \psi \right)_\mu = \phi *^D_\mu \psi$
\cite [Theorem 4.2] {BBEW96a}.

\section {Invariance properties and reduction}
\label {InvSec}

In \cite [p. 368]{BBEW96a} the notion of a `quantum moment map' was 
introduced for this situation (see e.~g. \cite {Xu96} for a more general 
discussion) and it was shown that $S_DJ = D(\lambda/x)J$ is a quantum moment
map for the $U(1)$-action for the star product $*^D$, i.~e. 
$D(\lambda/x)J$ induces the same group action on the quantum level as $J$ 
does on the classical level, i.~e. for all $F \in \mathcal A$ 
we have
\BEQ {QMM}
    F *^D S_DJ - S_DJ *^D F = \frac{\im\lambda}{2} \{ F, J \} .
\EEQ    
The classical observable algebra of the reduced system $C^\infty (\CP)$ 
can be thought as the quotient of the $U(1)$-invariant functions 
on $\C$ by those which vanish on the `constraint surface' 
$J^{-1} (\{\mu\})$. The later ideal is just the ideal generated by $J-\mu$ 
and thus one might have the idea that the quantum version works 
analogously: indeed in \cite[Prop. 4.1] {BBEW96a} it was shown that this 
is the case for $D=1$. We shall now consider the case of arbitrary series 
$D$. Then we have the already mentioned two possibilities: We can take the 
(left) ideal generated by $S_DJ - \mu$ in which case we can simply apply the
equivalence transformation $S_D$ to prove completely analogously to the proof
in \cite {BBEW96a} that the quotient space is isomorphic to the
functions on $\CP$ with the star product already obtained for $D=1$.
More interesting is hence the other possibility: we define the left-ideal
\BEQ {Ideal}
    \mathcal J^D_\mu := \mathcal A^0 *^D 
    \left( S_DJ - D\left(\frac{\lambda}{-2\mu}\right) \mu\right)
    \subset \mathcal A^0
\EEQ    
generated by $S_DJ - D(\lambda/(-2\mu))\mu$ which is in fact a two-sided ideal
due to (\ref {QMM}). Note that we have deformed both the classical 
momentum map {\em and} the momentum value in order to define the same 
constraint surface as in the classical case.
We now shall describe the quotient $\mathcal A^0 / J^D_\mu$:
\begin {lemma}
Denote by $\mathcal B := C^\infty (\CP)[[\lambda]]$ the vector space of 
the reduced observables then for any $F \in \mathcal A^0$ we have
\begin {enumerate}
\item $F_\mu = 0$ iff $F \in \mathcal J^D_\mu$.
\item $\mathcal B \cong \mathcal A^0 / J^D_\mu$ with the isomorphism
      $\mathcal B \ni \phi \mapsto [\pi^*\phi] \in \mathcal A^0/J^D_\mu$
      and its inverse 
      $\mathcal A^0 / J^D_\mu \ni [F] \mapsto F_\mu \in \mathcal B$.
\end {enumerate}
\end {lemma}      
\begin {proof}
Clearly $F_\mu = 0$ iff there exists a $G_0 \in \mathcal A^0$ such that 
$F = (J-\mu)G_0$ by Hadamard's trick 
(see e.~g. \cite [p. 366]{BBEW96a}). Since $S_DJ = J + \ldots$ we have
$F - (S_DJ - D(\lambda/(-2\mu))\mu)G_0 = \lambda F_1$
with some $F_1 \in \mathcal A^0$ and clearly $(F_1)_\mu = 0$ since 
$(S_DJ)_\mu = D(\lambda/(-2\mu))\mu$. Thus we can prove inductively that
$F_\mu = 0$ iff there exists a 
$G = G_0 + \lambda G_1 + \cdots \in \mathcal A^0$ such that 
$F = (S_DJ - D(\lambda/(-2\mu))\mu)G$. But since $S_DJ$ is radial we can
replace the pointwise product by the $*^D$ product which proves the first 
part. The second part follows directly.
\end {proof}
\begin {corollary}
The linear isomorphism 
$\mathcal B \ni \phi \mapsto [\pi^* \phi] \in \mathcal A^0 / J^D_\mu$
is an algebra isomorphism if $\mathcal B$ is equipped with the star 
product $*^D_\mu$ as in (\ref {RedStarDef}) and 
$\mathcal A^0/\mathcal J^D_\mu$ with the usual quotient algebra structure.
\end {corollary}
Since the star product algebra $(\mathcal B, *^D_\mu)$ is isomorphic 
to the quotient $\mathcal A^0/J^D_\mu$ and since $J^D_\mu$ is the quantum 
analogue of the classical vanishing ideal of functions vanishing on the 
classical constraint surface one can indeed  speak of $*^D_\mu$ as a 
reduced star product coming form the star product $*^D$ on $\C$.

\section {Non-equivalence of the reduced star products}
\label {NonSec}

Since in $\C$ all the star products $*^D$ are equivalent this construction 
raises the question whether the reduced star products are still equivalent in 
the sense of equivalence between star products. As we shall see by 
comparing the star products $*^D_\mu$ for varying series $D$ this is not the 
case. Using the concrete formula (\ref {HomProd}) for the $*^D$-product 
of homogeneous functions we easily obtain the following lemma by 
direct computation:
\begin {lemma}
Let $D, D'$ be two formal power series starting with $1$ and denote by 
$C, C'$ the corresponding inverse power series. Assume that $c_r = c_r'$ 
for $r = 1, \ldots, k-1$ and $c_k \ne c_k'$. Then the bidifferential 
operators of the corresponding star products $*^D$ and $*^{D'}$ coincide 
$K^D_r = K^{D'}_r$ for $r = 1, \ldots, k$ and in order $k+1$ we have
\BEQ {DiffKk}
    K^D_{k+1} - K^{D'}_{k+1} = (c_k - c_k') M_1 .
\EEQ
\end {lemma}
\begin {corollary} \label {DiffCor}
Under the same preconditions as in the preceeding lemma we have for 
the bidifferential operators in $*^D_\mu$ and $*^{D'}_\mu$
\BEQ {DiffTKr}
    \wt K^D_r = \wt K^{D'}_r  \mbox { for } r \le k
    \quad
    \mbox { and }
    \quad
    \wt K^D_{k+1} - \wt K^{D'}_{k+1} = (c_k - c'_k) \wt M_1 . 
\EEQ
\end {corollary}

Now it is known (see e.~g. \cite [Prop 3.7]{BCG97}) that if two 
{\em equivalent} star products on a 
symplectic manifold coincide up to order $k$ then the antisymmetric part
of their difference in order $k+1$ is a one-differential operator which 
can be written as $\Omega (X_f, X_g)$ with an {\em exact} two-form $\Omega$ 
where $X_f, X_g$ denote the Hamiltonian vector fields of the functions. 
But since in our case the antisymmetric part of $\wt M_1$ is just $\im/2$ 
times the Poisson bracket on $\CP$ which corresponds to the {\em non-exact} 
Fubini-Study form as two-form we immidiately have the following theorem:
\begin {theorem}
Let $D,D'$ be two complex formal power series starting with $1$ and let 
$*^D_\mu$ and $*^{D'}_\mu$ be the corresponding star product on $\CP$ 
according to (\ref {RedStarDef}). Then $*^D_\mu$ is equivalent to 
$*^{D'}_\mu$ iff $D=D'$ and any star product on $\CP$ is equivalent to some 
$*^D_\mu$. 
\end {theorem}
\begin {proof}
The non-equivalence for $D \ne D'$ follows easily form corollary 
\ref {DiffCor}, the fact that the Fubini-Study form is not exact and 
\cite [Prop 3.7] {BCG97}. Since $H^2 (\CP)$ is 
one-dimensional the possible equivalence classes are parametrised 
by $\lambda H^2 (\CP)[[\lambda]]$ (for fixed Poisson bracket, 
see e.~g. \cite {NT95a,BCG97}) which is clearly in bijection to the 
antisymmetric one-differential part determined by $c_kM_1$ since $c_k$ can 
be chosen arbitrarily (relative to the reference star product with $D=1$) 
in order $\lambda^{k+1}$.
\end {proof}

Remarks: In our example the choice $D=1$ with the classical value $\mu$ 
is clearly prefered since only for this choice the corresponding 
star product $*^D$ is strong $U(1)$-invariant, i.~e. the quantum momentum 
map coincides with the classical one. Nevertheless one can think of more 
general situation where no strong invariant star products are available.
Here one might weaken the strong invariance to the existence 
of a quantum momentum map and hence there might occur some subtilities in 
the choice of the star product and the choice of the quantum momentum 
value leading to non-equivalent star products for the quotient. 
Hence it would be very interesting to examine these aspects of reduction 
in more general situations. A good starting point for this programme 
should perhaps be Fedosov's reduction scheme \cite {Fed97}

\section* {Acknowledgements}

We would like to thank Martin Bordemann, Moshe Flato, Joachim Schirmer
and Daniel Sternheimer for a motivating discussion and
Martin Bordemann for carefully reading the manuscipt 
and useful comments.

\end {document}